\documentclass[a4paper,12pt]{article}
\usepackage{graphicx,amsmath,latexsym,amssymb,amsthm,cite,enumerate}
\usepackage[top=0.8in, bottom=1in, left=1.2in, right=1in]{geometry}

\newtheorem{thm}{Theorem}[section] 
\newtheorem{lem}[thm]{Lemma}
\newtheorem{cor}[thm]{Corollary}
\newtheorem{prop}[thm]{Proposition}
\theoremstyle{definition}
\newtheorem{defi}[thm]{Definition}
\newtheorem{ex}[thm]{Example}
\newtheorem{rem}[thm]{Remark}

\theoremstyle{plain}

\newcommand{\bR}{{\mathbb{R}}}

\begin{document}            

\def\au{Wanchai Tapanyo\textsuperscript{1}, Wachiraphong Ratiphaphongthon\textsuperscript{2},\\ Areerat Arunchai\textsuperscript{3}}

\def\inst{\normalsize
	\textsuperscript{1,3}Division of Mathematics and Statistics,  Nakhon Sawan Rajabhat University, Nakhon Sawan, Thailand,\\
	\textsuperscript{2}Department of Mathematics, Naresuan University, Phisanulok, Thailand\\
}

\def\ml{\textsuperscript{1}wanchai.t@nsru.ac.th, \textsuperscript{2}wachirapongr58@nu.ac.th,
\textsuperscript{3}areerat.a@nsru.ac.th}

\def\ti{On Completion of \textit{C}*-algebra-valued metric spaces}

\def\abs{
    The concept of a $ C ${\normalfont *}-algebra-valued metric space was introduced in 2014. It is a generalization of a metric space by replacing the set of real numbers by a $ C ${\normalfont *}-algebra. In this paper, we show that $ C ${\normalfont *}-algebra-valued metric spaces are cone metric spaces in some point of view which is useful to extend results of the cone case to $ C ${\normalfont *}-algebra-valued metric spaces. Then the completion theorem of $ C ${\normalfont *}-algebra-valued metric spaces is obtained. Moreover, the completion theorem of $ C ${\normalfont *}-algebra-valued normed spaces is verified and the connection with Hilbert $ C ${\normalfont *}-modules, generalized inner product spaces, is also provided.
}

\def\kw{$ C ${\normalfont *}-algebra-valued metric space, Cone metric space, Normed space, Inner product space, Completion}
	
\begin{center}
	~\vskip 2em
	{\LARGE\bf\ti}
	\vskip 1.5em
	{\large\au}
	\vskip 0.5em
	{\inst}
	\vskip 0.5em
	{Email: \ml}
	\vskip 3em
	{\Large\sc Abstract}\vskip 1em
\end{center}\normalsize

\abs 


\section{Introduction}

A metric space is one of attractive objects in mathematics which plays an important role in various branches of mathematics.
It is a nonempty set $ X $ together with a distance function $ d:X\times X\to \mathbb{R} $, which is often called a metric on $ X $.
Plenty of research papers study various kinds of spaces generalized from the definition of a metric space in different directions.
Some authors remove or change initial properties of a metric space while others change the values of the distance function to be in generalized sets of real or complex numbers, such as, a Banach space or a $ C ${\normalfont *}-algebra which can be seen in \cite{ConeHUANG2007} and \cite{CStarMa2014}, respectively.

The concept of a $ C ${\normalfont *}-algebra-valued metric space was first introduced in 2014 by Z. Ma and others. For this space the distance function was replaced by a function valued in a $ C ${\normalfont *}-algebra $ \mathbb{A} $.  If we consider the set of all positive elements $ \mathbb{A}_+ $ of $ \mathbb{A} $ as a cone of $ \mathbb{A} $. A $ C ${\normalfont *}-algebra-valued metric space is, in fact, a cone metric space which was introduced in 2004 by L.-G. Huang and others, see more details about a cone metric space in \cite{ConeHUANG2007}.

The main purpose of this research is to study the completion for $ C ${\normalfont *}-algebra-valued metric spaces and a $ C ${\normalfont *}-algebra-valued normed spaces. We verify some facts and use them to extend the results from others in \cite{ConeCompleteAbdeljawad2010}. Then we discuss relationships between $ C ${\normalfont *}-algebra-valued metric spaces and Hilbert $ C ${\normalfont *}-modules, generalized inner product spaces whose scalar fields are replaced by some $ C ${\normalfont *}-algebras.

The rest of the paper is organized as follows. 
In section 2 we derive the important definitions and theorems used to obtain our results. 
In section 3 We discuss on $ C ${\normalfont *}-algebra-valued metric and normed spaces and the relation to cone metric spaces. In section 4 the connection to Hilbert $ C ${\normalfont *}-modules is provided.

\section{Preliminaries}
This section provides a brief review of basic knowledge used in this research which can be found in \cite{CStarMa2014,CStarMurphy1990,ConeAnsari2017,ConeCompleteAbdeljawad2010,ConeHUANG2007}. We start with the definition of $ C ${\normalfont *}-algebras and some necessary related properties. Then we mention the definition of $ C ${\normalfont *}-algebra-valued and cone metric spaces and some previous results provided in other research papers.

\begin{defi}
	An \emph{algebra} is a vector space $\mathbb{A}$ together with a bilinear map $\mathbb{A}^2 \to \mathbb{A},$ $(a,b) \mapsto ab$, such that $a(bc)=(ab)c$ for all $a,b,c \in \mathbb{A}$.
\end{defi}

\begin{defi}
	An \emph{involution} on an algebra $\mathbb{A}$ is a conjugate-linear map $a \mapsto a^{*}$ on $\mathbb{A}$, such that $a^{**} = a$ and $(ab)^{*}=b^{*}a^{*}$ for all $a,b \in \mathbb{A}$. The pair $(\mathbb{A},*)$ is called an \emph{involutive algebra} or \emph{$*$-algebra}.
\end{defi}

\begin{defi}
	A \emph{Banach $\ast$-algebra} is a $\ast$-algebra $\mathbb{A}$ together with a complete norm such that $ \|ab\|\leq \|a\|\|b\| $ and $\lVert a^{*} \rVert = \lVert a \rVert$ for every $a, b \in \mathbb{A}$.
\end{defi}

\begin{defi}
	A \emph{$ C ${\normalfont *}-algebra} is a Banach $\ast$-algebra such that $\lVert a^{*}a \rVert = \lVert a \rVert^{2} $ for every $a \in \mathbb{A}$. If $ \mathbb{A} $ admits a unit $ I $ ($ aI= Ia= a $ for every $ a\in \mathbb{A} $) such that $ \|I\|= 1 $, we call $ \mathbb{A} $ a \emph{unital} \emph{$ C ${\normalfont *}-algebra}.
\end{defi}

\begin{rem}
	The word ``unital'' is also used for other terminologies, for example, unital algebra, unital $ * $-algebra, unital Banach $ * $-algebra. If a norm is not defined, the condition $ \|I\|= 1 $ will be omitted.
\end{rem}

\begin{defi}
	A \emph{homomorphism} is a linear map $ f $ from an algebra $ \mathbb{A} $ to an algebra $ \mathbb{B} $ such that $ f(ab)= f(a)f(b) $. An \emph{isomorphism} is a bijective homomorphism. In the case that $ \mathbb{A} $ and $ \mathbb{B} $ are $ * $-algebra, a \emph{$ * $-homomorphism (resp. $ * $-isomorphism)} is a homomorphism (resp. isomorphism) $ f: \mathbb{A}\to\mathbb{B} $ preserving adjoints, that is, $ f(a^*)= f(a)^* $ for every $ a\in \mathbb{A} $. In addition, if norms are defined on $ \mathbb{A} $ and $ \mathbb{B} $, the word ``isometric'' will be added before all of the terminologies to indicate that $ \|f(a)\|= \|a\| $.
\end{defi}

Consider the Cartesian $ (\mathbb{A},\mathbb{C}) $ of $ \mathbb{A} $ and the complex plane $ \mathbb{C} $, it is $ * $-algebra together with a unit $ (0,1) $ under componentwise operations for addition and involution, and the multiplication defined by
\[ (a,\alpha)(b,\beta)= (ab+\beta a+ \alpha b,\alpha\beta), \]
for every $ a,b\in\mathbb{A} $ and every $ \alpha,\beta\in\mathbb{C} $. By the norm obtained in \cite[Theorem 2.1.6]{CStarMurphy1990} we obtain that $ (\mathbb{A},\mathbb{C}) $ becomes a unital $ C ${\normalfont *}-algebra. The $ C ${\normalfont *}-algebra $ \mathbb{A} $ can be embedded in a unital $ C ${\normalfont *}-algebra $ (\mathbb{A},\mathbb{C}) $ by the injective isometric $ * $-homomorphism defined by 
	$$ a\mapsto (a,0). $$ 
The Cartesian $ (\mathbb{A},\mathbb{C}) $ is called the \emph{unitization} of $ \mathbb{A} $ and denoted by $ \widetilde{\mathbb{A}} $. Therefore, we may consider $ \mathbb{A} $ as a $ C ${\normalfont *}-subalgebra of $ \widetilde{\mathbb{A}} $.

For an element $ a $ of a unital algebra $ \mathbb{A} $, we say that $ a $ is \emph{invertible} if there is an element $ b\in \mathbb{A} $ such that $ ab= I= ba $. We denote by $ \mathrm{Inv}(\mathbb{A}) $ the set of all invertible elements of $ \mathbb{A} $. We define the \emph{spectrum} of $ a $ to be the set $ \sigma(a)=\sigma_\mathbb{A}(a) \{\lambda\in\mathbb{C}: \lambda I- a\notin\mathrm{Inv}(A) \} $. If $ \mathbb{A} $ is nonunital, we define $ \sigma_\mathbb{A}(a)= \sigma_{\widetilde{\mathbb{A}}}(a) $. The following definitions are about characterization of elements in $ C$*-algebra.

\begin{defi}
	An element $a$ of a $ * $-algebra $\mathbb{A}$ is called \emph{self-adjoint} or \emph{hermitian} if $a^{*} = a$.
	The set of all hermitian elements of $\mathbb{A}$ is denoted by $\mathbb{A}_h$. If $ \mathbb{A} $ is a $ C ${\normalfont *}-algebra, a self-adjoint element $ a\in\mathbb{A} $ with $ \sigma(a)\subseteq [0,+\infty) $ is called \emph{positive} and the set of all positive elements of $ \mathbb{A} $ is denoted by $ \mathbb{A}_+ $.
\end{defi}

If $ \mathbb{A} $ is a $ C ${\normalfont *}-algebra, $ \mathbb{A}_h $ becomes a partially ordered set by defining $ a \leq b $ to mean $ b-a\in \mathbb{A}_+ $. It is obvious that $ 0_\mathbb{A}\leq a $ precisely for $ a\in\mathbb{A}_+ $ where $ 0_\mathbb{A} $ is a zero in $ \mathbb{A} $. Thus, we may write $ 0_\mathbb{A}\leq a $ to indicate that $ a $ is positive.

\begin{prop}[Murphy] \label{Prop:Decomposition}
	Let $ \mathbb{A} $ be a $ C ${\normalfont *}-algebra. Then for each $ x\in \mathbb{A} $ there is a unique pair of hermitian elements $ a, b\in \mathbb{A} $ such that $ x= a+ bi $. More precisely, $ a= \frac{1}{2}(x+x^*) $ and $ b= \frac{1}{2i}(x-x^*)  $.
\end{prop}

\begin{thm}\cite[Theorem 2.2.1]{CStarMurphy1990}\label{Thm:root}
	Let $ a $ be a positive element of a $ C ${\normalfont *}-algebra $ \mathbb{A} $. Then there is a unique $ b\in \mathbb{A}_+ $ such that $ b^2= a $.
\end{thm}

By the previous theorem we can define the square root of the positive element $ a $ to be the element $ b $, we denote it by $ a^{1/2} $. The theorem below is a brief review of some necessary properties for positive elements of a $ C ${\normalfont *}-algebra, see more details in \cite{CStarMurphy1990}.

\begin{lem}\label{Lem:SumOfPositive}
	The sum of two positive elements in a $ C ${\normalfont *}-algebra are positive.
\end{lem}

\begin{thm}\label{thm:PropertiesOfPositive}
	Let $ \mathbb{A} $ be a $ C ${\normalfont *}-algebra. The the following properties are satisfied.
	\begin{enumerate}
		\item Suppose that $ \mathbb{A} $ is unital and $ a\in\mathbb{A} $ is hermitian. If $ \|a-tI\|\leq t $ for some $ t\in \mathbb{R} $, then $ a $ in positive. In the reverse direction, for every $ t\in\mathbb{R} $, if $ \|a\|\leq t $ and $ a $ is positive, then $ \|a-tI\|\leq t $
		\item For every $ a, b, c\in \mathbb{A}_h $, $ a\leq b $ implies $ a+c\leq b+c $,
		\item For every real numbers $ \alpha, \beta\geq 0 $ and every $ a, b\in \mathbb{A}_+ $, $ \alpha a +\beta b\in \mathbb{A}_+ $,
		\item $ A_+= \{a^*a : a\in\mathbb{A}\} $,
		\item If $ a,b\in A_h $ and $ c\in A $, then $ a\leq b $ implies $ c^*ac\leq c^*bc $,
		\item If $ 0_\mathbb{A}\leq a\leq b $, then $ \|a\|\leq \|b\| $.
	\end{enumerate}
\end{thm}

\begin{lem} \label{Lem:root2}
	Let $ \gamma=\alpha + \beta i \in\mathbb{C} $ and $ a\in\mathbb{A}_+ $. Then $ ((\alpha^2 + \beta^2)a)^{1/2}= |\gamma|a^{1/2} $.
	\begin{proof}
		It is obvious that $ |\gamma|a^{1/2} $ is positive. Consider 
		$$ (|\gamma|a^{1/2})^2= |\gamma|^2(a^{1/2})^2= (\alpha^2 +\beta^2)a. $$
		By Theorem \ref{Thm:root} , we have $ ((\alpha^2 + \beta^2)a)^{1/2}= |\gamma|a^{1/2} $.
	\end{proof}
\end{lem}

\begin{thm} \label{Thm:RootPresOder}
	Let $ a,b\in \mathbb{A}_+ $. Then $ a\leq b $ implies $ a^{1/2}\leq b^{1/2} $.
\end{thm}

Next, we provide the definitions of a $ C ${\normalfont *}-algebra-valued metric space, convergent sequences and Cauchy sequences in the space which are our main study.

\begin{defi}
	Let $ X $ be a nonempty set and $ d:X\times X\to \mathbb{A} $ be a function satisfying the following properties:
\begin{enumerate}[\indent (C1)]
	\item $ d(x,y)\geq 0 $,
	\item $ d(x,y)=0 $ if and only if $ x=y $,
	\item $ d(x,y)= d(y,x) $,
	\item $ d(x,y)\leq d(x,z)+d(z,y) $,
\end{enumerate}
for every $ x, y, z\in X $. We call the function $ d $ a $ C ${\normalfont *}-algebra-valued metric and call the triple $ (X,\mathbb{A},d) $ a $ C ${\normalfont *}-algebra-valued metric space. 
\end{defi}

We know that every $ C ${\normalfont *}-algebra $ \mathbb{A} $ can be embedded in $ \widetilde{\mathbb{A}} $ which is a unital $ C ${\normalfont *}-algebra. This means that we can work on $ \widetilde{\mathbb{A}} $ instead. In other words, an $ \mathbb{A} $-valued metric $ d $ of the space $ (X,\mathbb{A},d) $ is an $ \widetilde{\mathbb{A}} $-valued metric as concluded in the remark below.

\begin{rem}\label{rem:UnitisationValued}
	A $ C ${\normalfont *}-algebra-valued metric space $ (X,\mathbb{A},d) $ is a $ C ${\normalfont *}-algebra-valued metric space $ (X,\widetilde{\mathbb{A}},d) $.
\end{rem}

Therefore, we will assume $ \mathbb{A} $ to be unital in our research. The following statements are definitions of convergent and Cauchy sequences in a $ C ${\normalfont *}-algebra-valued metric space which are defined in \cite[Definition 2.2]{CStarMa2014}. We change some inequality in the definitions to correspond them to other similar definitions we use frequently.

\begin{defi}\label{ConvergentCauchyC*}
	Let $ (X,\mathbb{A},d) $ be a $ C ${\normalfont *}-algebra-valued metric space. A sequence $ \{x_n\} $ in $ X $ is said to \emph{converge} to an element $ x\in X $ (with respect to $ \mathbb{A} $) if and only if for every $ \varepsilon> 0 $ there is a positive integer $ N $ such that for every integer $ n\geq N $ we have $ \|d(x_n,x)\|<\varepsilon $. In this case we write $ \lim_{n\to\infty} x_n= x $, and say that the sequence $ \{x_n\} $ is \emph{convergent}.
	
	A sequence $ \{x_n\} $ in $ X $ is said to be \emph{Cauchy} (with respect to $ \mathbb{A} $) If and only if for every $ \varepsilon > 0 $ there is a positive integer $ N $ such that for every integer $ n,m\geq N $ we have $ \| d(x_n,x_m) \| < \varepsilon $.

We say that a $ C^\ast $-algebra-valued metric space $ (X, \mathbb{A},d) $ is \emph{complete} if every Cauchy sequence (with respect to $ \mathbb{A} $) is convergent.
\end{defi}

Next, we discuss cone metric spaces. We start with a cone of a real Banach Space which was introduced in \cite{ConeHUANG2007}. The definition is different one from \cite{ConeAnsari2017} which allows a cone to be trivial.

\begin{defi}
	Let $ \mathbb{E} $ be ba real Banach space. A nonempty closed subset $ P $ of $ \mathbb{E} $ is called a cone if and only if it satisfies the following properties:
\begin{enumerate}[\indent(P1)]
	\item $ P\neq\{0\} $,
	\item For every real numbers $ \alpha, \beta\geq 0 $ and every $ a, b\in P $, $ \alpha a +\beta b\in P $,
	\item If $ x\in P $ and $ -x\in P $, then $ x= 0 $.
\end{enumerate}
\end{defi}

Now we can define a partial order $ \leq $ on $ \mathbb{E} $ with respect to $ P $ by $ x\leq y $ to mean $ y-x\in P $. We write $ x< y $ to indicate that $ x\leq y $ and $ x\neq y $, and write $ x\ll y $ if $ y-x\in \mathrm{Int}(P) $.

A cone $ P $ is said to be \emph{normal} if and only if there exists a positive real number $ K $ such that for every $ x, y\in \mathbb{E} $, $ 0\leq x\leq y $ implies $ \|x\|\leq K\|y\| $. The following proposition is a consequence of Theorem \ref{thm:PropertiesOfPositive}. $ \mathbb{A}_+ $ is a cone in the sense of the preceding definition.

\begin{prop}
	$ \mathbb{A}_+ $ is a cone of a unital $ C ${\normalfont *}-algebra $ \mathbb{A} $.
	\begin{proof}
		The proof is immediate from Lemma \ref{Lem:SumOfPositive} and Theorem \ref{thm:PropertiesOfPositive}. Let $ \{x_n\} $ be a sequence in $ \mathbb{A}_+ $ converging to $ x\in\mathbb{A} $. Since $ \mathbb{A}_h $ is closed in $ \mathbb{A} $ and $ \mathbb{A}_+\subseteq\mathbb{A_h} $, we have $ x\in\mathbb{A}_h $. To show that $ \mathbb{A}_+ $ is closed we need to show that $ x\in\mathbb{A}_+ $. 
		
		Since $ \{x_n\} $ is convergent, it is certainly bounded. Then there is a positive real number $ t $ such that $ \|x_n\|\leq t $ for every $ n\in\mathbb{N} $. We know that $ x_n $ is positive for every $ n\in\mathbb{N} $. Thus, Theorem \ref{thm:PropertiesOfPositive} implies that $ \|x_n-tI\|\leq t $ for every $ n\in\mathbb{N} $. Consider
			\[ \|x-tI\| \leq \|x_n-x\|+\|x_n-tI\|\leq \|x_n-x\|+t. \]
		This implies that $ \|x-tI\|\leq t $. Since $ x $ is hermitian, again by Theorem \ref{thm:PropertiesOfPositive} we have $ x\in\mathbb{A}_+ $. Therefore, $ \mathbb{A}_+ $ is closed in $ \mathbb{A} $.  
	\end{proof}
\end{prop}

\begin{defi}
	Let $ X $ be a nonempty set and $ d:X\times X\to \mathbb{E} $ be a function satisfying the following properties:
	\begin{enumerate}[\indent (M1)]
		\item $ d(x,y)\geq 0 $,
		\item $ d(x,y)=0 $ if and only if $ x=y $,
		\item $ d(x,y)= d(y,x) $,
		\item $ d(x,y)\leq d(x,z)+d(z,y) $,
	\end{enumerate}
	for every $ x, y, z\in X $. We call the function $ d $ a \emph{cone metric} and call the pair $ (X,d) $ a \emph{cone metric space}. 
\end{defi}

Consider a $ C ${\normalfont *}-algebra $ \mathbb{A} $. If the scalar filed is restricted to the set of real numbers, $ \mathbb{A} $ becomes a real Banach space. Thus, a $ C ${\normalfont *}-algebra-valued metric space becomes a cone metric space.


\begin{defi}\label{ConvergentCauchyCone}
	Let $ (X,d) $ be a cone metric space. A sequence $ \{x_n\} $ in $ X $ is said to \emph{converge} to  $ x\in X $ (with respect to $ \mathbb{E} $) if and only if for every $ c\in \mathbb{E} $ with $ c\gg 0 $ there is a positive integer $ N $ such that for every integer $ n\geq N $ we have $ d(x_n,x)\ll c $. In this case we write $ \lim_{n\to\infty} x_n= x $, and say that the sequence $ \{x_n\} $ is \emph{convergent}.
	
	A sequence $ \{x_n\} $ in $ X $ is said to be \emph{Cauchy} (with respect to $ \mathbb{E} $) If and only if for every $ c\in \mathbb{E} $ with $ c\gg 0 $ there is a positive integer $ N $ such that for every integer $ n,m\geq N $ we have $ d(x_n,x_m) \ll c $.
	
	We say that a cone metric space $ (X,d) $ is complete if every Cauchy sequence (with respect to $ \mathbb{E} $) is convergent.
\end{defi}

\begin{lem}\label{lemEQUIV}
	Let $ (X,d) $ be a cone metric space together with a normal cone. A sequence $ \{x_n\} $ converses to $ x\in X $ if and only if $ \lim_{n\to+\infty}d(x_n,x)= 0 $. A sequence $ \{x_n\} $ is Cauchy if and only if $ \lim_{n,m\to+\infty}d(x_n,x_m)= 0 $. 
\end{lem}

\begin{defi}
	Let $ X $ be a vector space over the real field and $ \|\cdot\|: X\to \mathbb{E} $ be a function. A pair $ (X, \|\cdot\|) $ is called a \emph{cone normed space} if $ \|\cdot\| $ satisfies the following properties:
	\begin{enumerate}
		\item $ \|x\|= 0_\mathbb{E} $ if and only if $ x= 0_X $,
		\item $ \|\alpha x\|= |\alpha|\|x\| $,
		\item $ \|x + y\|\leq \|x\|+ \|y\| $,
	\end{enumerate}
	for every $ x, y\in X $ and every scalar $ \alpha $.
\end{defi}

\begin{thm}\label{thm:ConeCompletion}
	Let $ (X,d) $ be a cone metric space over a normal cone. Then there is a complete cone metric space $ (X^s,d^s) $ which has a dense subspace $W$ isometric to $X$. The space $X^s$ is unique except for isometries.
\end{thm}

\begin{thm}\label{thm:ConeNormCompletion}
	Let $ (X,\|\cdot\|) $ be a cone normed space over a normal cone. Then there is a cone Banach space $(X^s,\|\cdot\|^s) $ which has a dense subspace $W$ isometrically isomorphic to $X$. The space $X^s$ is unique except for isometric isomorphisms.
\end{thm}

The two results above are completion theorems obtained in \cite{ConeCompleteAbdeljawad2010}. We apply the first one to obtain our results. The isometry mentioned in the first theorem is a bijective mapping $ T:X\to Y $ between cone metric spaces preserving distances, that is,
\[ d_X(x,y)= d_Y(Tx,Ty), \]
for every $ x,y\in X $, where $ d_X $ and $ d_Y $ are metrics on $ X $ and $ Y $, respectively. Properties of the mapping $ T $ are different from those of the ordinary version only the values of $ d $ and $ d^s $ which are not real numbers. The second theorem is rewritten from the original version of cone normed spaces. The word ``isomorphism'' refers to a bijective linear operator between cone normed spaces and the word ``isometric'' indicates that the isomorphism is a cone-norm-preserving. In \cite{ConeCompleteAbdeljawad2010}, an isomorphism between cone normed spaces is always cone-norm-preserving.

Concepts of isometries and of $ C ${\normalfont *}-algebra-values metric spaces and isometric isomorphisms of $ C ${\normalfont *}-algebra-values normed spaces will be provided in the next section with more general than those of the cone version.


\section{Completion of \textit{C}*-algebra-valued metric\\ and normed spaces}

In this section we verify that a $ C ${\normalfont *}-algebra-valued metric space can be embedded in a complete $ C ${\normalfont *}-algebra-valued metric space as a dense subspace. The theorem in a version of a $ C ${\normalfont *}-algebra-valued normed space is also provided. We apply the fact that the $ C ${\normalfont *}-algebra-valued metric (resp. normed) spaces are cone metric (resp. normed) spaces to extend the results from the previous results in \cite{ConeCompleteAbdeljawad2010}. To work with a cone metric space, we need to assume that the interior of a cone is nonempty. However, this property does not generally occur for a $ C ${\normalfont *}-algebra as we show in the series of examples below.

\begin{ex}
	Let $ \mathbb{A} $ be a complex plane $ \mathbb{C} $. 
	Then $ \mathbb{A}_+= [0,\infty) $, so $ \mathrm{Int}(\mathbb{A}_+) $ is empty in $ \mathbb{C} $. Observe that $ \mathrm{Int}(\mathbb{A}_+) $ is not empty in $ \mathbb{R} $, the set of self-adjoint elements of $ \mathbb{C} $.\qed
\end{ex}

\begin{ex}
	In this example we consider $ \mathbb{A} $ as a $ C ${\normalfont *}-algebra of all bounded complex sequences $ \ell^\infty $ with the operators defined as follows:
	\begin{align*}
	(\xi_n)+ (\eta_n)	
	&= (\xi_n+ \eta_n)\\
	(\xi_n)(\eta_n)
	&= (\xi_n\eta_n)\\
	\lambda(\xi_n)
	&= (\lambda \xi_n)\\
	(\xi_n)^*
	&= (\mathbb{A}r{\xi}_n)\\
	\|(\xi_n)\|
	&= \sup_{n\in\mathbb{N}} |\xi_n|
	\end{align*}
	for every $ (\xi_n), (\eta_n)\in \ell^\infty $ and every $ \lambda\in \mathbb{C} $. We will show that $ \mathrm{Int}(\ell^\infty_+)=\emptyset $
	
	By the definitions of self-adjoint and positive elements of a $ C ${\normalfont *}-algebra, we have 
	\[ \ell^\infty_h
	= \big\{a\in\ell^\infty : a^*= a \big\}
	= \big\{(\xi_n)\in\ell^\infty : \xi_n\in\mathbb{R} \text{ for all } n\in\mathbb{N} \big\} \]
	and
	\[ \ell^\infty_+
	= \big\{a\in\ell^\infty_h : \sigma(a)\subseteq \mathbb{R}_+ \big\}
	= \big\{(\xi_n)\in\ell^\infty : \xi_n\in\mathbb{R}_+ \text{ for all } n\in\mathbb{N} \big\}. \]
	To show that $ \mathrm{Int}(\ell^\infty_+)=\emptyset $, we let $ a= (\xi_n)\in \ell^\infty_+ $ and $ \varepsilon> 0 $. Then choose $ b= (\xi_1 - i\frac{\varepsilon}{2}, \xi_2, \xi_3,\ldots) $. Clearly, $ b $ is in $ \ell^\infty\setminus\ell^\infty_+ $ such that $ \|a-b\|= \frac{\varepsilon}{2}< \varepsilon $. This implies that $ b\in B(a,\varepsilon) $, the open ball in $ \ell^\infty $ of radius $ \varepsilon $ centered at $ a $. Since $ \varepsilon $ is arbitrary, the element $ a $ is not an interior point of $ \ell^\infty_+ $. This situation occurs for every element of $ \ell^\infty_+ $, so we have $ \mathrm{Int}(\ell^\infty_+)= \emptyset $.\qed
\end{ex}

\begin{ex}[A $ C ${\normalfont *}-algebra-valued metric space with the empty interior of $ \mathbb{A}_+ $]\
	
	\indent	In this example we replace $ X $ and $ \mathbb{A} $ by $ \mathbb{C} $ and $ \mathbb{C}^2 $, respectively. By the same operators in the previous example, the space $ \mathbb{C}^2 $ can be considered as a $ C ${\normalfont *}-subalgebra of $ \ell^\infty $ with $ \mathrm{Int}(\mathbb{C}^2_+)=\emptyset $. Let $ d:\mathbb{C}\times\mathbb{C}\to \mathbb{C}^2 $ be a function defined by
	\[ d(a,b)= (|a-b|,\alpha|a-b|) \]
	such that $ \alpha>0 $ for every $ a, b\in\mathbb{C} $. Therefore, $ (\mathbb{C}, \mathbb{C}^2, d) $ is a $ C ${\normalfont *}-algebra-valued metric space.\qed
\end{ex}

Although the situation in the previous example can occur, we does not assume that the interior of $ \mathbb{A}_+ $ is not empty. This is the result of that there exists a suitable real Banach subspace of $ \mathbb{A} $ containing $ \mathbb{A}_+ $ with a nonempty interior under the topology on the Banach subspace restricted from $ \mathbb{A} $, and so, we will work on the subspace instead. We conclude this useful fact in the two following propositions.

\begin{prop}\label{prop:AhBanach}
	$ \mathbb{A}_h $ is a real Banach subspace of a $ C ${\normalfont *}-algebra $ \mathbb{A} $.
	\begin{proof}
		Since $\mathbb{A}_h \subseteq \mathbb{A}$, $0_\mathbb{A} \in \mathbb{A}_h$ and $(\alpha a + b)^* = \alpha a + b$ for all $\alpha \in \bR$ and $a,b \in \mathbb{A}_h$, we obtain that $\mathbb{A}_h$ is a real normed space. The completeness of $ \mathbb{A}_h $ can be obtained by verifying that $ \mathbb{A}_h $ is closed in $ \mathbb{A} $. Let $\{ a_n \}$ be a sequence in $\mathbb{A}_h$ converging to $ a\in \mathbb{A} $. Since $\lVert a_n - a \rVert = \lVert (a_n - a)^* \rVert = \lVert a_n^* - a^* \rVert = \lVert a_n - a^* \rVert$, we obtain that $a_n $ converges to $ a^*$. By the uniqueness of a limit of a convergent sequence, we have $a=a^*$, i.e. $a \in \mathbb{A}_h$. Therefore, $\mathbb{A}_h$ is closed in $\mathbb{A}$, and so $\mathbb{A}_h$ is a real Banach subspace of $ \mathbb{A} $.
	\end{proof}
\end{prop}

\begin{prop}\label{prop:NonemptyInterior}
	If $ \mathbb{A} $ is a unital $ C ${\normalfont *}-algebra, then $ \mathrm{Int}_{A_h}(\mathbb{A}_+)\neq\emptyset $. 
	\begin{proof}
		Let $ I $ be a unit of $ \mathbb{A} $ and $ B(I,1)= \{a\in\mathbb{A}_h:\|a-I\|<1\} $. Then Theorem \ref{thm:PropertiesOfPositive} implies that $ B(I,1)\subseteq \mathbb{A}_+ $. Hence, $ I\in \mathrm{Int}_{A_h}(\mathbb{A}_+) $, so $ \mathrm{Int}_{A_h}(\mathbb{A}_+)\neq\emptyset $.
	\end{proof}
\end{prop}

\begin{cor}
	If $\mathbb{A}$ is a unital $ C ${\normalfont *}-algebra and $\mathbb{A} = \mathbb{A}_h$, then $\mathrm{Int}(\mathbb{A}) \neq \emptyset$.
\end{cor}

Since $ \mathbb{A_+} $ is a cone of $ \mathbb{A} $, it is also a cone of $ \mathbb{A}_h $. We have known by Theorem \ref{thm:PropertiesOfPositive} that $ \|a\|\leq \|b\| $ for every positive elements $ a $ and $ b $ satisfying $ a\leq b $. Hence, $ \mathbb{A_+} $ is a normal cone of $ \mathbb{A}_h $. This is a fact occurring in every $ C^\ast $-algebra.  Therefore, a $ C^\ast $-algebra-valued metric space $ (X,\mathbb{A},d) $ is a cone metric space $ (X,\widetilde{\mathbb{A}}_h,d) $ with a normal cone $ \widetilde{\mathbb{A}}_+ $ such that $ \mathrm{Int}_{\widetilde{\mathbb{A}}_h}(\widetilde{\mathbb{A}}_+)\neq\emptyset $. Finally, we obtain Lemma \ref{lemEQUIV} in a version of a $ C^\ast $-algebra-valued metric space $ (X,\mathbb{A},d) $, equivalent definitions of convergent and Cauchy sequences, stated in the following theorem.

\begin{thm}\label{thmEquiConvCstar}
	Let $ (x_n) $ be a sequence in a $ C ${\normalfont *}-algebra-valued metric space $ (X,\mathbb{A},d) $. Then the following statements are satisfied.
	\begin{enumerate}
		\item $ (x_n) $ converges to $ x\in\mathbb{X} $ (in the sense of Definition \ref{ConvergentCauchyC*}) if and only if for every $ c\in \widetilde{\mathbb{A}}_h $ with $ c\gg 0 $ there is a positive integer $ N $ such that for every integer $ n\geq N $ we have $ d(x_n,x)\ll c $.
		\item $ (X_n) $ is Cauchy (in the sense of Definition \ref{ConvergentCauchyC*}) if and only if for every $ c\in \widetilde{\mathbb{A}}_h $ with $ c\gg 0 $ there is a positive integer $ N $ such that for every integer $ n,m\geq N $ we have $ d(x_n,x_m)\ll c $.
	\end{enumerate}
\begin{proof}
	We prove only the convergent case, the other can be proved similarly. Suppose that $ (x_n) $ converges to an element $ x $ of $ (X,\mathbb{A},d) $. Then $ (x_n) $ converges to an element $ x $ of $ (X,\widetilde{\mathbb{A}},d) $, and so, converges in $ (X,\widetilde{\mathbb{A}}_h,d) $. Then the forward implication is obtained after applying Lemma \ref{lemEQUIV}. For the converse implication, we suppose that the condition holds. Then Lemma \ref{lemEQUIV} implies that $ \lim_{n\to+\infty}\|d(x_n,x)\|_{\widetilde{\mathbb{A}}_h}= 0 $. Since $ d(x_n,x) $ belongs to $ \mathbb{A} $, we have $ \lim_{n\to+\infty}\|d(x_n,x)\|_{{\mathbb{A}}}= 0 $. Therefore, $ (x_n) $ converges to an element $ x $ of $ (X,\mathbb{A},d) $.
\end{proof}
\end{thm}

Before verifying the completion theorem, we need to introduce some necessary definitions first. Let $ (X, \mathbb{A}, d_X) $ and $ (Y, \mathbb{B}, d_Y) $ be $ C ${\normalfont *}-algebra-valued metric spaces. A mapping $ T: X\to Y $ is call an \emph{isometry} if there exists an isometric $ * $-isomorphism $ f: \mathbb{A}\to \mathbb{B} $ such that
\[ f(d_X(x, y))= d_Y(T(x), T(y)), \]
for every $ x, y\in X $. The space $ (X, \mathbb{A}, d_X) $ and $ (Y, \mathbb{B}, d_Y) $ are said to be \emph{isometric} if there exists a bijective isometry from $ (X, \mathbb{A}, d_X) $ to $ (Y, \mathbb{B}, d_Y) $.

\begin{prop}
	An isometry between $ C ${\normalfont *}-algebra-valued metric spaces is always injective.
	\begin{proof}
		Suppose that $ (X, \mathbb{A}, d_X) $ and $ (Y, \mathbb{B}, d_Y) $ are $ C ${\normalfont *}-algebra-valued metric spaces and $ T $ is an isometry from $ X $ to $ Y $. Without loss of generality, we may assume that $ \mathbb{B}= \mathbb{A} $. Let $ x,y\in X $ such that $ T(x)=T(y) $. Then $ d_X(x, y)= d_Y(T(x), T(y))= 0_\mathbb{A} $, so $ x=y $. Therefore, $ T $ is injective.
	\end{proof}
\end{prop}

The definition of denseness of a subset of a topological space is determined using neighborhoods. In the case of a $ C ${\normalfont *}-algebra-valued metric space, we provide an equivalent definition using only open balls in the space. 

\begin{defi}\label{defClosure}
	Let $(X,\mathbb{A},d)$ be a C*-algebra-valued metric space and $M$ be a subset of $X$. For any $ \varepsilon> 0 $, we define 
		$$ B(x, \varepsilon)= \{y\in X: \|d(x,y)\|_\mathbb{A}< \varepsilon\} $$
	Let $ M $ be a subset of $ X $, the set of all limit points or \emph{closure} of $M$ is determined by 
		$$\mathrm{Cl}(M) = \{ x \in X : B(x, \varepsilon) \cap M \neq \emptyset \text{ for every } \varepsilon > 0\}.$$
	If $ \mathrm{Cl}(M)= X $, we say that $ M $ is \emph{dense} in $ X $.
\end{defi}

Because of Theorem \ref{thmEquiConvCstar}, an equivalent definition of closure of the set $ M $ is obtained, that is,
	\[ \mathrm{Cl}(M) = \{ x \in X : B_1(x, c) \cap M \neq \emptyset \text{ for every } c \gg 0\}, \]
where $ B_1(x,c)= \{y\in X: d(x,y)< c\} $ 
with $ c\in\mathbb{A} $ such that $ c\gg 0 $.


We have shown that any $ C ${\normalfont *}-algebra-valued metric space $ (X, \mathbb{A}, d) $ can be considered as the cone metric space $ (X, \widetilde{\mathbb{A}}_h, d) $ with the normal cone $ \widetilde{\mathbb{A}}_+ $ such that $ \mathrm{Int}_{\widetilde{\mathbb{A}}_h}(\widetilde{\mathbb{A}}_+)\neq\emptyset $. Thus, we can work on the cone metric space instead, and obtain the completion of $ (X, \mathbb{A}_h, d) $ after applying Theorem \ref{thm:ConeCompletion}. Since the values of $ d $ belong to $ \mathbb{A} $, the $ C ${\normalfont *}-algebra-valued metric space $ (X, \mathbb{A}, d) $ is actuary contained in the acquired space as a dense subspace. We conclude this result in the following theorem.

\begin{thm}[Completion of $ C ${\normalfont *}-algebra-valued metric spaces] \label{thm:MetricCompletion}\ \\	
	For any $ C ${\normalfont *}-algebra-valued metric space $(X, \mathbb{A}, d)$, there exists a complete $ C ${\normalfont *}-algebra-valued metric space $(X^s, \mathbb{A}, d^s)$ which contains a dense subspace $W$ isometric with $X$. The space $X^s$ is unique except for isometries.
	\begin{proof}
		We know that $(X, \mathbb{A}, d)$ is a cone metric space $(X, \widetilde{\mathbb{A}}_h, d)$ containing the normal cone $ \widetilde{\mathbb{A}}_+ $ such that $ \mathrm{Int}_{\widetilde{\mathbb{A}}_h}(\widetilde{\mathbb{A}}_+)\neq\emptyset $. Then Theorem \ref{thm:ConeCompletion} implies that there is a complete cone metric space $ (X^s, \widetilde{\mathbb{A}}_h, {d}^s) $ which contains a dense subspace $ W $ isometric with $ X $. We see that $ (X^s,\widetilde{\mathbb{A}}, {d}^s) $ is also a $ C ${\normalfont *}-algebra-valued metric space. We will verify that $ d^s $ is an $ \mathbb{A} $-valued metric for $ X^s $, in fact, after taking the composition with the inverse of the mapping $ a\mapsto (a,0) $ from $ \mathbb{A} $ to $ \widetilde{\mathbb{A}} $.
		
		Let $ x,y\in X^s $. Since $ W $ is dense in $ X^s $, there exist sequences $ \{x_n\} $ and $ \{y_n\} $ in $ W $ converging to $ x $ and $ y $, respectively. By continuity of $ {d}^s $, we have 
				\[ {d}^s(x,y)= \lim_{n\to \infty} {d}^s(x_n,y_n). \]
		Let $ T $ be a bijective isometry of a cone metric space from $ W $ to $ (X, \widetilde{\mathbb{A}}_h, d) $. Then 
			\[ {d}^s(x_n,y_n)= d(T(x_n),T(y_n)) \in \mathbb{A}, \]
		for every $ n\in\mathbb{N} $. Since $ \mathbb{A} $ is closed in $ \widetilde{\mathbb{A}} $, we have 
			\[ d^s(x,y)= \lim_{n\to \infty} d^s(x_n,y_n) \in \mathbb{A}. \] 
		This implies that $ d^s $ is an $ \mathbb{A} $-valued metric for $ X^s $. Let $ (\widehat{X},\mathbb{B},\widehat{d}) $ be another $ C ${\normalfont *}-algebra-valued metric space which contains a dense subspace $ \widehat{W} $ isometric with $ X $. Then there is a bijective isometry $ \widehat{T} $ from $ X $ to $ \widehat{W} $. Thus, $ \widehat{T}\circ T $ is a bijective isometry from $ W $ to $ \widehat{W} $. Therefore, $ \widehat{T}\circ T $ can be extended to be a bijective isometry from $ X^s $ to $ \widehat{X} $ after applying denseness of $ W $ and $ \widehat{W} $ in $ X^s $ and $ \widehat{X} $, respectively. This means that the space $ X^s $ exists uniquely except for isometries. 
	\end{proof}
\end{thm}

Next, we focus on a $ C ${\normalfont *}-algebra-valued normed space. We know that any incomplete normed space is embeddable in another complete normed space. In \cite{ConeCompleteAbdeljawad2010}, the concept of completion is also extended to the case of a cone normed space. The author defined a cone normed space and verified the existence of its completion. Therefore, the result in the case of $ C ${\normalfont *}-algebra-valued normed spaces is obtained directly from the sake of Proposition \ref{prop:AhBanach}  and Proposition \ref{prop:NonemptyInterior} similar to the case of a $ C ${\normalfont *}-algebra-valued metric space.

Let $ X $ be a vector space over the real or complex fields and $ \mathbb{A} $ be a $ C ${\normalfont *}-algebra. A triple $ (X, \mathbb{A}, \|\cdot\|) $ is called a $ C ${\normalfont *}-algebra-valued normed space if $ \|\cdot\| $ is a function from $ X $ to $ A_+ $ satisfying the following properties:
\begin{enumerate}
	\item $ \|x\|= 0_\mathbb{A} $ if and only if $ x= 0_X $,
	\item $ \|\alpha x\|= |\alpha|\|x\| $,
	\item $ \|x + y\|\leq \|x\|+ \|y\| $,
\end{enumerate}
for every $ x, y\in X $ and every scalar $ \alpha $. Notice that $ 0_\mathbb{A} $ and $ 0_X $ are zeros in $ \mathbb{A} $ and $ X $ respectively.

By the definition of a $ C ${\normalfont *}-algebra-valued norm, we can investigate that the function $ d: X\times X\to \mathbb{A} $ determined by $ d(x,y)= \|x-y\| $ is a $ C ${\normalfont *}-algebra-valued metric. We call it the \emph{$ C ${\normalfont *}-algebra-valued metric} induced by the norm $ \|\cdot\| $ We conclude this fact in the proposition below

\begin{prop}\label{prop:MetricIndByNorm}
	A $ C ${\normalfont *}-algebra-valued normed space $ (X, \mathbb{A}, \|\cdot\|) $ is a $ C ${\normalfont *}-algebra-valued metric space with a metric $ d: X\times X\to \mathbb{A} $ given by $ d(x,y)= \|x-y\| $.
\end{prop}

A complete $ C ${\normalfont *}-algebra-valued normed space under the metric defined above is called a \emph{$ C ${\normalfont *}-algebra-valued Banach Space}. In the next example, we show that every commutative $ C ${\normalfont *}-algebra is a $ C ${\normalfont *}-algebra-valued normed space. We provide a lemma before verify this claim.

\begin{lem} \label{Lem:root3}
	Let $ A $ be commutative $ C ${\normalfont *}-algebra. Then $ \mathbb{A}_h $ is a closed $ * $-subalgebra of $ \mathbb{A} $ over the real field. Moreover, if $ a,b\in \mathbb{A}_+ $, then $ ab\in \mathbb{A}_+ $ and $ (ab)^{1/2}= a^{1/2}b^{1/2} $.
	\begin{proof}
		Since $ \mathbb{A} $ is commutative, $ (ab)^*= a^*b^*= ab $ for every $ a,b\in\mathbb{A}_h $. Combine with Proposition \ref{prop:AhBanach}, $ \mathbb{A}_h $ is a real $ * $-subalgebra of $ \mathbb{A} $. Now, suppose that $ a,b\in\mathbb{A}_+ $ Theorem \ref{thm:PropertiesOfPositive} implies that $ a= c^*c $ for some $ c\in \mathbb{A} $. Thus, we have $ 0_\mathbb{A}= c^*0_\mathbb{A}c\leq c^*bc= c^*cb= ab $, so $ ab $ is positive. By the same way, $ a^{1/2}b^{1/2} $ is also positive. Since $ (a^{1/2}b^{1/2})^2= ab $, Theorem \ref{Thm:root} implies that $ a^{1/2}b^{1/2}= (ab)^{1/2} $.
	\end{proof}
\end{lem}

\begin{ex}\label{ex:CStarValNormForA}
	Let $ \mathbb{A} $ be a commutative $ C ${\normalfont *}-algebra and $ X= \mathbb{A} $. We know from Proposition \ref{Prop:Decomposition} that every element $ x\in \mathbb{A} $ can be uniquely decomposed as $ x= a+ bi $ for some $ a, b\in \mathbb{A}_h $. Then we define $ \|\cdot\|_0: X \to \mathbb{A}_+ $ by 
	\[ \|x\|_0= (a^2 + b^2)^{1/2}. \]
	We will show that $ (X,\|\cdot\|_0,\mathbb{A}) $ is a $ C ${\normalfont *}-algebra-valued normed space.
	
	Since $ a $ and $ b $ are hermitian, Theorem \ref{thm:PropertiesOfPositive} implies that $ a^2 $ and $ b^2 $ are positive. Thus, $ (a^2+ b^2)^{1/2} $ is also positive after applying Lemma \ref{Lem:SumOfPositive} and Theorem \ref{Thm:root}, respectively. This shows that $ \|\cdot\|_0 $ is a well-defined function with its values in $ \mathbb{A}_+ $. Since $ x= 0_X $ if and only if $ a=b=0_X $, we obtain that $ \|x\|_0= 0_\mathbb{A} $ if and only if $ x= 0_X $. Next, let $ \gamma\in \mathbb{C} $. Then $ \gamma= \alpha+ \beta i $ where $ \alpha, \beta\in \mathbb{R} $. Hence, $ \gamma x= (\alpha +\beta i)(a+bi)= (\alpha a- \beta b) + (\beta a+ \alpha b)i $, so
	\begin{align*}
	\|\gamma x\|^2_0 &= (\alpha a- \beta b)^2 + (\beta a + \alpha b)^2\\
	&= \alpha^2a^2 + \beta^2b^2 + \beta^2a^2 + \alpha^2 b^2\\
	&= (\alpha^2 + \beta^2)(a^2 + b^2).
	\end{align*} 
	Theorem \ref{Thm:root} and Lemma \ref{Lem:root2} imply that $ \|\gamma x\|_0= ((\alpha^2 +\beta^2)(a^2+b^2))^{1/2}= |\alpha|\|x\|_0 $.
	
	Finally, we prove the triangle inequality. Let $ y\in X $ be uniquely represented by $ c+ di $ where $ c, d\in \mathbb{A}_h $. Consider 
	\begin{align*}
	\|x+y\|^2_0 &= \|(a+c)+(b+d)i\|^2_0\\
	&= (a+c)^2+(b+d)^2\\
	&= (a^2+ 2ac + c^2) + (b^2+ 2bd+ b^2)\\
	&= (a^2 + b^2 + c^2 + d^2) + 2(ac + bd),
	\end{align*}
	and
	\begin{align*}
	(\|x\|_0+\|y\|_0)^2 &= \|x\|^2_0 + 2\|x\|_0\|y\|_0 + \|y\|^2_0\\
	&= (a^2+ b^2) + 2(a^2+b^2)^{1/2}(c^2+d^2)^{1/2} + (c^2+ d^2)\\
	&= (a^2 + b^2 + c^2 + d^2) + 2(a^2+b^2)^{1/2}(c^2+d^2)^{1/2}.
	\end{align*}
	We obtain by Theorem \ref{Thm:RootPresOder} that 
	\begin{center}
		$ \|x+y\|_0\leq \|x\|_0+\|y\|_0 $ whenever $ \|x+y\|^2_0\leq (\|x\|_0+\|y\|_0)^2. $
	\end{center}
	Thus, to complete this proof, we need to show that $ ac + bd\leq (a^2+b^2)^{1/2}(c^2+d^2)^{1/2} $. 
	
	We may see that $ 0_\mathbb{A}\leq (ad- bc)^2= (ad)^2 - 2abcd +(bc)^2 $, so $ 2abcd\leq (ad)^2 + (bc)^2 $. Therefore,
	\begin{align*}
	(ac+bd)^2 &= (ac)^2 + 2abcd + (bd)^2\\
	&\leq (ac)^2 + (ad)^2 + (bc)^2 + (bd)^2\\
	&= (a^2+ b^2)(c^2+ d^2).
	\end{align*}
	Theorem \ref{Thm:RootPresOder} implies $ ((ac + bd)^2)^{1/2}\leq ((a^2+ b^2)(c^2+ d^2))^{1/2} $. Then apply Theorem \ref{Thm:root} and Lemma \ref{Lem:root3} to the left and right sides of the inequality, respectively, so we obtain $ ac + bd\leq (a^2+b^2)^{1/2}(c^2+d^2)^{1/2} $. Now, the proof of triangle inequality is completed. Consequently, $ \|\cdot\|_0 $ is an $ \mathbb{A} $-valued norm for $ \mathbb{A} $.
	\qed
\end{ex}

Consider a linear operator between normed spaces.  It is an isometry if and only if it is norm-preserving. Thus, the isometric property may be replaced by the norm-preserving property to define isometric isomorphisms of normed spaces. We do the same for $ C ${\normalfont *}-algebra-valued normed spaces. Let $ (X, \mathbb{A}, \|\cdot\|_X) $ and $ (Y, \mathbb{B}, \|\cdot\|_Y) $ be $ C ${\normalfont *}-algebra-valued normed spaces. A linear operator $ T: X\to Y $ is called an \emph{isometric} if there exists a $ * $-isomorphism $ f: \mathbb{A}\to \mathbb{B} $ such that
\[ f(\|x\|_X)= \|T(x)\|_Y, \]
for every $ x\in X $. A bijective isometric linear operator is called an \emph{isometric isomorphism}. We say that the spaces $ (X, \mathbb{A}, \|\cdot\|_X) $ and $ (Y, \mathbb{B}, \|\cdot\|_Y) $ are \emph{isometrically isomorphic} if there exists an isometric isomorphism from $ (X, \mathbb{A}, \|\cdot\|_X) $ to $ (Y, \mathbb{B}, \|\cdot\|_Y) $.

\begin{lem}\label{lem:InducedNorm}
	Let $ (W,\mathbb{A},d) $ and $ (X,\mathbb{A},\|\cdot\|_X) $ be $ C ${\normalfont *}-algebra-valued metric and normed spaces, respectively. Assume that $ T $ is a bijective isometry from $ W $ and $ X $. Then the following statements are satisfied.
	\begin{enumerate}
		\item $ (W,\mathbb{A}, \|\cdot\|_W) $ is a $ C ${\normalfont *}-algebra-valued normed space such that $ \|u\|_W= \|Tu\|_X $ for every $ u\in W $, so $ T $ is a norm-preserving operator. Moreover, $ d(u,v)= \|u-v\|_W $ for every $ u,v\in W $. 
		\item $ T $ becomes a linear operator from $ W $ to $ X $.
	\end{enumerate}
	\begin{proof}
		We define the additive operation $ \oplus $ and the scalar multiplication $ \odot $ on $ W $ by 
			$$ u\oplus v:= T^{-1}(Tu+Tv) \text{ and } \alpha\odot u:= T^{-1}(\alpha Tu),$$ 
		for every $ u,v\in W $ and every scalar $ \alpha\in\mathbb{C} $. We will show that the operations satisfy all axioms of a vector space. Let $ u, v, w\in W $ and $ \alpha,\beta\in\mathbb{C} $ be scalars.
		\begin{enumerate}[V1)]
			\item Closure property of $ \oplus $ and $ \odot $: Clearly, $ u\oplus v, \alpha\odot u\in W $.
			\item Associativity of $ \oplus $: $ \begin{aligned}[t]
				(u\oplus v)\oplus w &= T^{-1}\Big(T\big(T^{-1}(Tu+Tv)\big)+Tw\Big)\\
					&= T^{-1}(Tu+Tv+Tw)\\
					&= T^{-1}\Big(Tu+T\big(T^{-1}(Tu+Tv)\big)\Big)\\
					&= u\oplus (v\oplus w)
			\end{aligned} $
			
			\item Commutativity of $ \oplus $: $ \begin{aligned}[t]
				u\oplus v &= T^{-1}(Tu+Tv)\\
					&= T^{-1}(Tv+Tu)\\
					&= v\oplus u
			\end{aligned} $
			
			\item\label{item:identity} Let $ 0_X $ be the identity of $ X $. We show that $ T^{-1}(0_X) $ is the identity of $ W $ under $ \oplus $. For any element $ u\in W $, we have
			\begin{align*}
				T^{-1}(0_X)\oplus u &= u\oplus T^{-1}(0_X)\\
					&= T^{-1}\Big(Tu+T\big(T^{-1}(0_X)\big)\Big) \\
					&= T^{-1}(Tu+0_X)= u.
			\end{align*}

			\item We show that $ T^{-1}(-Tu) $ is the inverse of $ u $ in $ W $ under $ \oplus $,
			\begin{align*}
				T^{-1}(-Tu)\oplus u &= u\oplus T^{-1}(-Tu)\\
					&= T^{-1}\Big(Tu+ T\big(T^{-1}(-Tu)\big)\Big)\\
					&= T^{-1}(0_X).
			\end{align*}
					
			\item Compatibility of $\odot$ with multiplication of the filed $ \mathbb{C} $:
			\begin{align*}
				\alpha\odot(\beta\odot u) &= T^{-1}\Big(\alpha T\big(T^{-1}(\beta Tu)\big)\Big)\\
					&= T^{-1}\big(\alpha (\beta Tu)\big)\\
					&= T^{-1}\big((\alpha\beta) Tu\big)\\
					&= (\alpha\beta)\odot u.
			\end{align*}

			\item Distributivity of $ \odot $ with respect to $ \oplus $:
				\begin{align*}
					\alpha\odot(u\oplus v) &= T^{-1}\Big(\alpha T\big(T^{-1}(Tu+Tv)\big)\Big)\\
						&= T^{-1}(\alpha Tu+ \alpha Tv)\\
						&= T^{-1}\Big(T\big(T^{-1}(\alpha Tu)\big)+ T\big(T^{-1}(\alpha Tv)\big)\Big)\\
						&= T^{-1}(\alpha Tu)\oplus T^{-1}(\alpha Tv)\\
						&= (\alpha\odot u)\oplus (\alpha\odot v)
				\end{align*}
			
			\item Distributivity of $ \odot $ with respect to addition of the field $ \mathbb{C} $:
				\begin{align*}
					(\alpha+\beta)\odot u &= T^{-1}\big((\alpha+\beta)Tu\big)\\
						&= T^{-1}(\alpha Tu+ \beta Tu)\\
						&= T^{-1}\Big(T\big(T^{-1}(\alpha Tu)\big)+ T\big(T^{-1}(\beta Tu)\big)\Big)\\
						&= T^{-1}(\alpha Tu)\oplus T^{-1}(\beta Tu)\\
						&= (\alpha\odot u)\oplus (\beta\odot u).
				\end{align*}
			
			\item Identity under $ \odot $:
				$ 1\odot u= T^{-1}(1Tu)= u $.
		\end{enumerate}
		Now, $ (W,\oplus,\odot) $ is a vector space over the field $ \mathbb{C} $. Next, we verify that $ \|\cdot\|_W $ is an $ \mathbb{A} $-valued norm on $ W $. Let $ 0_W= T^{-1}(0_X) $, the identity of $ W $. Clearly, $ \|0_W\|= \|0_X\|_X= 0_\mathbb{A} $. In the reverse direction, we assume that $ \|u\|= 0_\mathbb{A} $. Then $ \|Tu\|_X= 0_\mathbb{A} $. Thus, $ Tu= 0_X $, so $ u= T^{-1}(0_X)= 0_W $. Moreover, we have
			 $$ \|\alpha\odot u\|= \|T\big(T^{-1}(\alpha Tu)\big)\|_X= |u|\|Tu\|_X= |\alpha|\|u\| $$
		and	
			 $$ \|u\oplus v\|= \|T^{-1}(Tu+Tv)\|= \|Tu+Tv\|_X\leq \|Tu\|_X+\|Tv\|_X= \|u\|+\|v\|. $$
	Therefore, $ (W,\mathbb{A},\|\cdot\|_W) $ is a $ C ${\normalfont *}-algebra-valued normed space such that
	\begin{align*}
		d_W(u,v) &= d_X(Tu,Tv)\\
			&= \|Tu-Tv\|_X\\
			&= \|T^{-1}(Tu-Tv)\|\\
			&= \|u\oplus T^{-1}(-Tv)\|\\
			&= \|u\oplus (-v)\|.
	\end{align*}
	In addition, for every $ u,v\in W $ and every $ \alpha,\beta\in\mathbb{C} $, we have
	\begin{align*}
		T\big((\alpha\odot u)\oplus v\big) &= T\Big(T^{-1}\big(T(\alpha\odot u)+T(\beta\odot v)\big)\Big)\\
			&= T(\alpha\odot u)+T(\beta\odot v)\\
			&= T\big(T^{-1}(\alpha Tu)\big)+ T\big(T^{-1}(\beta Tv)\big)\\
			&= \alpha Tu+ \beta Tv
	\end{align*}	
Therefore, $ T $ is a linear operator. The proof of the lemma is now complete.
	\end{proof}
\end{lem}

\begin{lem}\label{lem:BecomeBanach}
	Let $ (W,\mathbb{A},\|\cdot\|_W) $ be a $ C ${\normalfont *}-algebra-valued normed space contained as a subspace of a complete $ C ${\normalfont *}-algebra-valued metric space $ (X^s,\mathbb{A},d^s) $. Assume that $ W $ is dense in $ X^s $. Then $ X^s $ becomes $ C ${\normalfont *}-algebra-valued Banach space.
	\begin{proof}
		Let $ x^s, y^s\in X^s $ and $ \alpha\in\mathbb{C} $ be a scalar. Then there are sequences $ \{x_n\} $ and $ \{y_n\} $ in $ W $ converging to $ x^s $ and $ y^s $, respectively. Consider
			\[ \|(x_n+y_n)- (x_m+y_m)\|_W\leq \|x_n-x_m\|_W+\|y_n-y_m\|_W, \]
		and
			\[ \|\alpha x_n-\alpha x_m\|_W= |\alpha|\|x_n-x_m\|_W. \]
		These imply that $ \{x_n+y_n\} $ and $ \{\alpha x_n\} $ are Cauchy, so they converge in $ X^s $. We extend the addition and the scalar multiplication of $ W $ to $ X^s $ by $ x^s+y^s=\lim_{n\to\infty} (x_n+y_n) $ and $ \alpha x^s= \lim_{n\to\infty} \alpha x_n $. Next, we show that the extended operations are well-defined.  
		
		Assume that $ x'_n $ and $ y'_n $ are other sequences in $ W $ converging to $ x^s $ and $ y^s $, respectively. We define a sequence $ \{z_n\} $ by $ z_{2n-1}= x_n+y_n $ and $ z_{2n}= x'_n+y'_n $ for every $ n\in \mathbb{N} $. We see that
			\begin{align*}
				\|(x_n+y_n)-(x'_m+y'_m)\| &\leq \|x_n-x'_m\|_W+\|y_n-y'_m\|_W\\
						&\leq d_s(x_n,x^s)+d_s(x^s,x'_m)+d_s(y_n,y^s)+d_s(y^s,y'_m).
			\end{align*}
		Thus, $ \{z_n\} $ is Cauchy in $ X^s $, so it converges in $ X^s $. This implies that its subsequences $ \{z_{2n}\} $ and $ \{z_{2n-1}\} $ converge to the same limit. Therefore, the addition is well-defined. The proof for scalar multiplication can be done similarly. Compatibility of the operations with the axiom of a vector space can be obtained by considering sequences in $ W $. Finally, $ X^s $ becomes a vector space with the identity $ 0_W $ and the limit of a sequence $ \{-x_n\} $ as the inverse of $ x^s $. We see that $ \{\|x_n\|\} $ is a convergent sequence in $ \mathbb{A} $, so we put $ \|x^s\|= \lim_{n\to\infty} \|x_n\| $. Consequently, $ X^s $ is an $ \mathbb{A} $-valued normed space. 
	\end{proof}
\end{lem}

\begin{thm}[Completion of $ C ${\normalfont *}-algebra-valued normed spaces] \label{thm:ComCstarNorm}\ \\	
	For any $ C ${\normalfont *}-algebra-valued normed space $(X, \mathbb{A}, \|\cdot\|)$, there exists a $ C ${\normalfont *}-algebra-valued Banach space $(X^s, \mathbb{A}, \|\cdot\|^s)$ which contains a dense subspace $W$ isometric with $X$. The space $X^s$ is unique except for isometries.
	\begin{proof}
		Let $ d $ be the metric induced by the norm as in Proposition \ref{prop:MetricIndByNorm}. We obtain that the $ C ${\normalfont *}-algebra-valued normed space $ (X,\mathbb{A},\|\cdot\|) $ becomes an $ \mathbb{A} $-valued metric space.  Then apply Theorem \ref{thm:MetricCompletion} to obtain a complete $ C ${\normalfont *}-algebra-valued metric space $ (X^s,\mathbb{A},d^s) $ containing a dense metric subspace $ W $ isometric with $ X $. By using the previous two lemmas $ X^s $ becomes a $ C ${\normalfont *}-algebra-valued Banach space with an $ \mathbb{A} $-valued norm $ \|\cdot\|^s $ such that $ d^s(x^s,y^s)= \|x^s - y^s\|^s $ for every $ x^s, y^s\in X^s $.
	\end{proof}
\end{thm}

In Lemma \ref{lem:InducedNorm}, we show that the bijective isometry $ T $ from the space $(W, \mathbb{A}, d)$ to the space $(X, \mathbb{A}, \|\cdot\|_X)$ finally becomes a linear operator. Suppose that this situation occurs for other spaces $ \widehat{W} $ and $ \widehat{X}^s $ together with a bijective isometry $ \widehat{T}: \widehat{W}\to X $. Then the composition $ \widehat{T}^{-1}\circ T $ is a bijective linear operator from $ W $ to $ \widehat{W} $. By applying denseness of the spaces $ W $ and $ \widehat{W} $, we can extend the bijective linear operator to be an isometric isomorphism from $ X^s $ to $ \widehat{X}^s $. Now, we obtain another version of the preceding theorem stated in the corollary below by using isometric isomorphisms instead of isometries. The similar result studied in \cite{ConeCompleteAbdeljawad2010} is concluded in Theorem \ref{thm:ConeNormCompletion} for cone normed spaces.

\begin{cor}\label{cor:ComCstarNorm}
		For any $ C ${\normalfont *}-algebra-valued normed space $(X, \mathbb{A}, \|\cdot\|)$, there exists a $ C ${\normalfont *}-algebra-valued Banach space $ (X^s, \mathbb{A}, \|\cdot\|^s)$ which contains a dense subspace $W$ isometrically isomorphic with $X$. The space $X^s$ is unique except for isometric isomorphism.
\end{cor}

\section{Connection with Hilbert \textit{C}*-modules.}
This section provides certain relationships between concepts of a $ C ${\normalfont *}-algebra-valued metric space and an inner-product $ C ${\normalfont *}-module which is a generalization of an inner product space. The concept of inner-product $ C ${\normalfont *}-module was first introduced in \cite{KaplanskyIntroModule}, the study of I. Kaplansky in 1953, to develop the theory for commutative unital algebras. In the 1970s, the definition was extended to the case of noncommutative $ C ${\normalfont *}-algebra, see more details in \cite{PaschkeExtModule,RieffelExtModule}. Let $ \mathbb{A} $ be a $ C ${\normalfont *}-algebra and $ X $ be a complex vector space which is a right $ \mathbb{A} $-module with compatible scalar multiplication:
	\begin{equation}\label{eqn:CompatibleScalarPro}
		\alpha(xa)= (\alpha x)a= x(\alpha a),
	\end{equation}
for every $ \alpha\in\mathbb{C}, x\in X $ and $ a\in \mathbb{A} $.
The triple $ (X, \mathbb{A}, \langle\cdot,\cdot\rangle) $ is called an \emph{inner product $ \mathbb{A} $-module} if the mapping $ \langle\cdot,\cdot\rangle: X\times X \to \mathbb{A} $ satisfies the following conditions;
\begin{enumerate}
	\item $ \langle x, \alpha y+ \beta z\rangle= \alpha\langle x,y \rangle + \beta\langle x,z \rangle $,
	\item $\langle x,ya \rangle= \langle x,y \rangle a $,
	\item $ \langle y,x \rangle= \langle x,y \rangle^* $,
	\item $\langle x,x \rangle \geq 0_\mathbb{A}$,
	\item if $\langle x,x \rangle = 0_\mathbb{A}$, then $ x=0_X $,
\end{enumerate}
for every $ \alpha\in\mathbb{C} $ and every $ x, y\in X $. It is known that any inner product $ C ${\normalfont *}-module $ (X, \mathbb{A}, \langle\cdot,\cdot\rangle) $ is a norm space with a scalar-valued norm $ \|\cdot\|_m $ given by 
	$$ \|x\|_m= \|\langle x,x \rangle\|_\mathbb{A}^{1/2}, $$ 
for every $ x\in X $ where $ \|\cdot\|_\mathbb{A} $ is a norm on $ \mathbb{A} $. It is called a \emph{Hilbert $ C ${\normalfont *}-module} if the induced norm is complete.

Let $ \mathbb{A} $ be a commutative unital $ C ${\normalfont *}-algebra such that every nonzero element is invertible. We have $ \mathbb{A}= \mathbb{C}I $ where $ I $ is a unit of $ \mathbb{A} $. In this case an inner product $ C ${\normalfont *}-module is almost like a traditional inner product space that we can see easily. So the Cauchy-Schwarz inequality is also satisfied in context for a $ C ${\normalfont *}-algebra-valued inner product. In fact, we require only the values of the $ C ${\normalfont *}-algebra-valued inner product are invertible elements in the commutative unital $ C ${\normalfont *}-algebra. The inequality is proved in the following lemma.

\begin{lem}[Cauchy-Schwarz inequality]\ \\
	Let $ (X, \mathbb{A}, \langle\cdot,\cdot\rangle) $ be an inner product $ C ${\normalfont *}-module with a commutative unital $ C ${\normalfont *}-algebra $ \mathbb{A} $ such that every nonzero value of $ \langle \cdot, \cdot \rangle $ is invertible. Then
	\[ \langle x,y \rangle\langle y,x \rangle\leq \langle x,x \rangle\langle y,y \rangle, \]
	for every $ x,y\in X $. Moreover,
		\[ \|\langle x,y \rangle\|_0\leq \langle x,x \rangle^{1/2}\langle y,y \rangle^{1/2}, \]
	where $ \|\cdot\|_0 $ is the norm defined in Example \ref{ex:CStarValNormForA}.
	\begin{proof}
		Let $ x,y\in X $, $ a=\langle x,y \rangle $, $ b= \alpha I $ where $ \alpha\in\mathbb{R}_+ $ and $ I $ is a unit of $ \mathbb{A} $. Then
		\begin{align*}
		0_\mathbb{A} &\leq \langle xa- yb, xa- yb \rangle\\
		&= \langle xa,xa \rangle - \langle xa,yb \rangle - \langle yb,xa  \rangle + \langle yb,yb \rangle\\
		&= a^*\langle x,x \rangle a - a^*\langle x,y \rangle b - b^*\langle y,x \rangle a + b^*\langle y,y \rangle b\\
		&= \langle x,x \rangle a^*a - 2a^*ab + \langle y,y \rangle b^2.
		\end{align*}
		This implies that $ 2a^*ab\leq \langle x,x \rangle a^*a + \langle y,y \rangle b^2  $. 
		
		If $ \langle x,x \rangle= 0_\mathbb{A} $, then $ 2a^*a \leq \langle y,y \rangle b= \alpha \langle y,y \rangle $. This is true for every $ \alpha\in\mathbb{R}_+ $, so we have $ 2a^*a=0 $. Thus, $ \langle x,y \rangle\langle y,x \rangle\leq \langle x,x \rangle\langle y,y \rangle $. Now we assume that $ \langle x,x \rangle\neq 0_\mathbb{A} $. In this case, we let $ b= \langle x,x \rangle $. Hence, $ 2a^*ab\leq ba^*a + \langle y,y \rangle b^2  $, so $ a^*a\leq \langle y,y \rangle b  $. Therefore, $ \langle x,y \rangle\langle y,x \rangle\leq \langle x,x \rangle\langle y,y \rangle $.
		
		Next, we show that $ \|\langle x,y \rangle\|_0\leq \langle x,x \rangle^{1/2}\langle y,y \rangle^{1/2} $. By commutativity of $ \mathbb{A} $ and the representation 
			$$ \langle x,y \rangle= a+bi, $$ 
		for some $ a,b\in \mathbb{A}_h, $ we have	
		$ \|\langle x,y \rangle\|^2_0= \langle x,y \rangle\langle y,x \rangle\leq \langle x,x \rangle\langle y,y \rangle $. Then Theorem \ref{Thm:RootPresOder} implies that 
			$$ \|\langle x,y \rangle\|_0= (\|\langle x,y \rangle\|^2_0)^{1/2}\leq (\langle x,x \rangle\langle y,y \rangle)^{1/2}. $$ 
		Consider $ (\langle x,x \rangle^{1/2}\langle y,y \rangle^{1/2})^2= (\langle x,x \rangle^{1/2})^{2} (\langle y,y \rangle^{1/2})^2= \langle x,x \rangle\langle y,y \rangle. $
		Thus, 
			$$ (\langle x,x \rangle\langle y,y \rangle)^{1/2}= \langle x,x \rangle^{1/2}\langle y,y \rangle^{1/2}. $$ Therefore, 
			\begin{equation*}
			\|\langle x,y \rangle\|_0\leq \langle x,x \rangle^{1/2}\langle y,y \rangle^{1/2}. \qedhere
			\end{equation*}
	\end{proof}
\end{lem}

\begin{lem}
	Let $ a $ be a positive element of a $ C ${\normalfont *}-algebra $ \mathbb{A} $. Then $ \|a\|_0= a $ where $ \|\cdot\|_0 $ is the norm defined in Example \ref{ex:CStarValNormForA}
	\begin{proof}
		The definition of $ \|\cdot\|_0 $ implies that $ \|a\|_0= (a^2)^{1/2} $. Since $ a $ is the unique positive element such that $ a^2= a^2 $, we have $ (a^2)^{1/2}=a $. The proof is now completed.
	\end{proof}
\end{lem}

\begin{thm}
	Let $ (X, \mathbb{A}, \langle\cdot,\cdot\rangle) $ is an inner product $ C ${\normalfont *}-module. If $ \mathbb{A} $ is a commutative unital $ C ${\normalfont *}-algebra such that every nonzero value of $ \langle \cdot, \cdot \rangle $ is invertible, then $ X $ becomes a $ C ${\normalfont *}-algebra-valued normed space with a $ \mathbb{A} $-valued norm $ \|\cdot\| $ given by $ \|x\|= \langle x,x \rangle^{1/2} $, for every $ x\in X $.
	\begin{proof} Let $ x, y\in X $ and $ \alpha $ be a scalar. Since $ \langle x,x \rangle\in\mathbb{A}_+ $, $ \|x\|= \langle x,x \rangle^{1/2}\in\mathbb{A}_+ $.
		
		(1) Assume that $ \langle x,x \rangle^{1/2}=\|x\|= 0_\mathbb{A} $. Then $ \langle x,x \rangle = 0_\mathbb{A} $, so $ x= 0_X $. In reverse direction we assume that $ x= 0_X $. Then apply the second property in the definition of inner product $ C ${\normalfont *}-module and obtain $ \langle x,x \rangle= 0_\mathbb{A}= 0_\mathbb{A}^2 $, so $ \|x\|= \langle x,x \rangle^{1/2}= 0_\mathbb{A} $.
		
		(2) Consider $ \|\alpha x\|^2= \langle\alpha x, \alpha x \rangle= \alpha\bar{\alpha}\langle x,x  \rangle= |\alpha|^2\langle x,x  \rangle $, so $ \langle x,x \rangle= (\frac{1}{|\alpha|}\|\alpha x\|)^2 $. This implies that $ \|x\|= \langle x,x \rangle^{1/2}= \frac{1}{|\alpha|}\|\alpha x\| $, so $ |\alpha|\|x\|= \|\alpha x\| $.
		
		(3) Let us consider $ \| x-y \|^2 $. Clearly,
		\begin{align*}
		\| x-y \|^2
		&= \langle x-y, x-y \rangle\\
		&= \langle x,x \rangle - \langle x,y \rangle - \langle y,x \rangle + \langle y,y \rangle\\
		&= \|x\|^2 - \langle x,y \rangle - \langle y,x \rangle + \|y\|^2.
		\end{align*}
		Then take the norm $ \|\cdot\|_0 $ determined in Example \ref{ex:CStarValNormForA} to both sides of the equation and apply the two preceding lemmas, so we have
		\begin{align*}
		\| x-y \|^2
		&\leq \|x\|^2 + \|\langle x,y \rangle\|_0+ \|\langle y,x \rangle\|_0+ \|y\|^2\\
		&\leq \|x\|^2 + 2\|x\|\|y\|+ \|y\|^2\\
		&= (\|x\|+\|y\|)^2.
		\end{align*}
		This means that $ \| x-y \|= (\| x-y \|^2)^{1/2}\leq ((\|x\|+\|y\|)^2)^{1/2}= \|x\|+\|y\| $.\\
		From (1) to (3), we obtain that $ (X, \mathbb{A}, \|\cdot\|) $ is a $ C ${\normalfont *}-algebra-valued normed space.
	\end{proof}
\end{thm}

\begin{cor}
	Let $ (X, \mathbb{A}, \langle\cdot,\cdot\rangle) $ is an inner product $ C ${\normalfont *}-module with a commutative unital $ C ${\normalfont *}-algebra $ \mathbb{A} $ such that every nonzero element is invertible. Then $ X $ becomes a $ C ${\normalfont *}-algebra-valued normed space with the same norm determined in the preceding theorem.
\end{cor}

 In case the inner product $ C ${\normalfont *}-module $ (X, \mathbb{A}, \langle\cdot,\cdot\rangle) $ is a $ C ${\normalfont *}-algebra-valued normed space, and so a $ C ${\normalfont *}-algebra-valued metric space, we can consider whether the space is complete by using a $ C ${\normalfont *}-algebra-valued metric. The following theorem shows that these two definitions of completeness are identical in this situation.

\begin{thm}
	Assume that an inner product $ C ${\normalfont *}-module $ (X, \mathbb{A}, \langle\cdot,\cdot\rangle) $ is a $ C ${\normalfont *}-algebra-valued norm space with an $ \mathbb{A} $-valued norm $ \|\cdot\|_X $ induced by $ \langle\cdot,\cdot\rangle $. Then it is a Hilbert $ C ${\normalfont *}-module if and only if it is a $ C ${\normalfont *}-algebra-valued Banach space.
	\begin{proof}
		Let $ x $ be any element of $ X $ and $ \|\cdot\|_X $ be an $ \mathbb{A} $-valued norm on $ X $ induced by $ \langle\cdot,\cdot\rangle $. Since $ \|x\|_X^2=  \langle x,x \rangle $, we have
			
			\[ \big\|\langle x,x \rangle\big\|_\mathbb{A}= \big\|\|x\|_X^2\big\|_\mathbb{A}= \big\|\|x\|_X\big\|_\mathbb{A}^2. \]
		Thus,
			\[ \big\|\|x\|_X\big\|_\mathbb{A}= \big\|\langle x,x \rangle\big\|_\mathbb{A}^{1/2} =\|x\|_m. \]
		Then by Definition \ref{ConvergentCauchyC*} we obtain that the two concepts of convergence of any sequence $ \{x_n\} $ in $ X $ by $ \|\cdot\|_X $ and $ \|\cdot\| $ are equivalence. Therefore, $ X $ is a Hilbert $ C ${\normalfont *}-module if and only if it is a $ C ${\normalfont *}-algebra-valued Banach space.
	\end{proof}
\end{thm} 

The concept of completion is also extended to inner product $ C ${\normalfont *}-module. It is mentioned in \cite{HilbertModulesChristopher1995} that for any inner product $ C ${\normalfont *}-module $ X $ over a $ C ${\normalfont *}-algebra $ \mathbb{A} $, one can form its completion $ X^s $, a Hilbert $ \mathbb{A} $-module, using a similar way to the case of the scalar-valued inner product space. An $ \mathbb{A} $-valued inner product on $ X^s $ is constructed from one of $ X $ using the completeness of $ \mathbb{A} $.

By using the scalar-valued norm $ \|\cdot\|_m $, Corollary \ref{cor:ComCstarNorm} implies that there is a Banach space $ X^s $ which contains a dense subspace $ W $ isometrically isomorphic with $ (X,\|\cdot\|_m) $. Let $ T $ be an isometric isomorphism from $ W $ to $ X $. By the same argument of Lemma \ref{lem:InducedNorm}, the algebraic operation of a right $ \mathbb{A} $-module compatible with scalar multiplication \eqref{eqn:CompatibleScalarPro} on $ W $ can be induced by that on $ X $ through the mapping $ T $, that is,
	\[ ua := T^{-1}\big((Tu)a\big), \]
for every $ u\in W $ and every $ a\in \mathbb{A} $. An $ \mathbb{A} $-valued inner product for $ W $ can be induced in a similar way, that is,
	\[ \langle u,v \rangle_W= \langle Tu,Tv \rangle_X, \]
for every $ u,v\in W $. Now $ W $ becomes an inner product $ \mathbb{A} $-module and $ T $ is an isomorphism between inner product $ \mathbb{A} $-modules. Finally, we extend all the induced operations on $ W $ to $ X^s $ by the similar argument used in Lemma \ref{lem:BecomeBanach}. Let $ x^s, y^s\in X^s $. Then there exist a sequence $ \{x_n\} $ and $ \{y_n\} $ in $ W $ converging by the norm $ \|\cdot\|_m $ to $ x^s $ and $ y^s $, respectively. For every $ a\in\mathbb{A} $, we define an $ \mathbb{A} $-module operation and an $ \mathbb{A} $-valued inner product on $ X^s $ by
	\[ x^sa := \lim_{n\to\infty} x_na \qquad\text{and}\qquad \langle x^s,y^s\rangle:= \lim_{n\to\infty}\langle x_n,y_n\rangle. \]
It is the fact that norm $ \|\cdot\|_m $ makes $ W $ be a right normed $ \mathbb{A} $-module, that is, 
	$$ \|ua\|_m\leq \|u\|_m\|a\|_\mathbb{A}, $$ 
for every $ u\in W $ and every $ a\in\mathbb{A} $. Moreover, It follows from \cite[Proposition 1.1]{HilbertModulesChristopher1995} that 
	$$ \|\langle u,v\rangle\|_\mathbb{A}\leq \|u\|_m\|v\|_m $$
for every $ u, v\in W $. Therefore, the limits above exist and the operations for $ X^s $ is well-defined. Thus, $ X^s $ becomes a Hilbert $ \mathbb{A} $-module. We conclude this result in the remark below.

\begin{rem}
	The completion for any inner product $ C ${\normalfont *}-module $ (X,\mathbb{A},\langle\cdot,\cdot\rangle) $ exists. That is, there is a Hilbert $ C ${\normalfont *}-module $ (X^s,\mathbb{A},\langle\cdot,\cdot\rangle^s) $ containing $ W $ as an inner product $ C ${\normalfont *}-submodule such that $ W $ is isomorphic to $ X $ as $ \mathbb{A} $-valued inner product spaces.
\end{rem}

Assume that an inner product $ C ${\normalfont *}-module $ (X, \mathbb{A}, \langle\cdot,\cdot\rangle) $ is a $ C ${\normalfont *}-algebra-valued norm space with an $ \mathbb{A} $-valued norm $ \|\cdot\|_X $ induced by $ \langle\cdot,\cdot\rangle $. Then $ W $ is also a $ C ${\normalfont *}-algebra-valued norm space. In this case the completion of $ X $ can be constructed by using the norm $ \|\cdot\|_W $ instead of $ \|\cdot\|_m $ on $ W $. This is a result of the following identity,
	\[ \big\|\|u\|_W\big\|_\mathbb{A}= \|u\|_m, \]
for every $ u\in W $. Consequently, we also obtain the completion of an inner product $ C ${\normalfont *}-module $ X $ by applying the completion theorem for $ \mathbb{A} $-valued normed spaces if the $ \mathbb{A} $-valued norm induced by $ \langle\cdot,\cdot\rangle $ exists. An $ \mathbb{A} $-valued Inner product for the complete space can be induced from $ \langle\cdot,\cdot\rangle $ by using its continuity together with the Cauchy-Schwarz inequality. The continuity of $ \langle\cdot,\cdot\rangle $ is proved in the following theorem by applying the original version of the Cauchy-Schwarz inequality for $ X $, that is,
	\[ \langle y,x\rangle\langle x,y\rangle\leq \|\langle x,x\rangle\|_\mathbb{A}\langle y,y\rangle, \]
for every $ x,y\in X $. Since $ \langle y,x\rangle\langle x,y\rangle $ is positive, we have
	 $$ \|\langle y,x\rangle\|^2_\mathbb{A}=\|\langle y,x\rangle\langle x,y\rangle\|_\mathbb{A}\leq \|\langle x,x\rangle\|_\mathbb{A}\|\langle y,y\rangle\|_\mathbb{A}. $$
Therefore, the inequality below holds,
	\begin{equation}\label{CauchySchMo}
	\|\langle y,x\rangle\|_\mathbb{A}\leq \|\langle x,x\rangle\|_\mathbb{A}^{1/2}\|\langle y,y\rangle\|_\mathbb{A}^{1/2}.
	\end{equation}

\begin{thm}
	Let $ (X, \mathbb{A}, \langle\cdot,\cdot\rangle) $ be an inner product $ C ${\normalfont *}-module. Assume that the $ \mathbb{A} $-valued norm $ \|\cdot\| $ induced by $ \langle\cdot,\cdot\rangle $ exists. If $ x_n\to x $ and $ y_n\to y $ by the norm, then $ \langle x_n, y_n\rangle\to \langle x, y\rangle $ in $ \mathbb{A} $.
	\begin{proof}
		We apply the triangle inequality for $ \|\cdot\|_\mathbb{A} $ and then the inequality \eqref{CauchySchMo}, finally, rewrite the inner product in the form of the norm $ \|\cdot\| $, so we obtain
		\begin{align*}
			\|\langle x_n,y_n\rangle- \langle x,y\rangle\|_\mathbb{A} 
				&= \|\langle x_n,y_n\rangle-\langle x_n,y\rangle+ \langle x_n,y\rangle -\langle x,y\rangle\|_\mathbb{A}\\
				&\leq \|\langle x_n,y_n\rangle-\langle x_n,y\rangle\|_\mathbb{A} + \|\langle x_n,y\rangle -\langle x,y\rangle\|_\mathbb{A}\\
				&\leq \|\langle x_n,y_n-y\rangle\|_\mathbb{A} + \|\langle x_n-x,y\rangle\|_\mathbb{A}\\
				&\leq \|\langle x_n,y_n-y\rangle\|_\mathbb{A} + \|\langle x_n-x,y\rangle\|_\mathbb{A}\\
				&\leq \|\langle x_n,x_n\rangle\|^{1/2}_\mathbb{A} \|\langle y_n-y,y_n-y\rangle\|^{1/2}_\mathbb{A}\\
				&\quad +\|\langle x_n-x,x_n-x\rangle\|^{1/2}_\mathbb{A} \|\langle y,y\rangle\|^{1/2}_\mathbb{A}\\
				&\leq \left(\big\| \|x_n\|^2 \big\|^{1/2}_\mathbb{A}\right) \left(\big\| \| y_n-y\|^2 \big\|^{1/2}_\mathbb{A}\right)\\
				&\quad +\left(\big\| \|x_n-x\|^2 \big\|^{1/2}_\mathbb{A}\right) \left(\big\| \| y\|^2 \big\|^{1/2}_\mathbb{A}\right)\\
				&\leq \left( \big\| \|x_n\| \big\|_\mathbb{A}\right) \left(\big\| \| y_n-y\| \big\|_\mathbb{A}\right)  +\left( \big\| \|x_n-x\| \big\|_\mathbb{A}\right)  \left( \big\| \| y\| \big\|_\mathbb{A}\right).
		\end{align*}
	It not difficult to see that $ \big\| \|x_n\| \big\|_\mathbb{A} $ is bounded. Since $ x_n\to x $ and $ y_n\to y $ by the $ \mathbb{A} $-valued norm $ \|\cdot\| $, $ \|\langle x_n,y_n\rangle- \langle x,y\rangle\|_\mathbb{A}\to 0 $. Therefore, $ \langle x_n, y_n\rangle\to \langle x, y\rangle $ in $ \mathbb{A} $. 
	\end{proof}
\end{thm}
\begin{cor}
	Let $ (X, \mathbb{A}, \langle\cdot,\cdot\rangle) $ be an inner product $ C ${\normalfont *}-module with a commutative unital $ C ${\normalfont *}-algebra $ \mathbb{A} $ such that every nonzero value of $ \langle \cdot, \cdot \rangle $ is invertible. If $ x_n\to x $ and $ y_n\to y $ by the norm, then $ \langle x_n, y_n\rangle\to \langle x, y\rangle $ in $ \mathbb{A} $.
\end{cor}


\begin{thebibliography}{1}
	
	\bibitem{ConeHUANG2007}
	L.-G. Huang and X.~Zhang, ``Cone metric spaces and fixed point theorems of
	contractive mappings,'' {\em Journal of Mathematical Analysis and
		Applications}, vol.~332, no.~2, pp.~1468 -- 1476, 2007.
	
	\bibitem{CStarMa2014}
	Z.~Ma, L.~Jiang, and H.~Sun, ``$\uppercase{C}$*-algebra-valued metric spaces
	and related fixed point theorems,'' {\em Fixed Point Theory and
		Applications}, vol.~2014, p.~206, Oct 2014.
	
	\bibitem{ConeCompleteAbdeljawad2010}
	T.~Abdeljawad, ``Completion of cone metric spaces,'' {\em Hacettepe Journal of
		Mathematics and Statistics}, vol.~39, no.~1, pp.~67 -- 74, 2010.
	
	\bibitem{CStarMurphy1990}
	G.~J. Murphy, ``\uppercase{C}hapter 1 - elementary spectral theory and
	\uppercase{C}hapter 2 - $\uppercase{C}$*-algebras and hilbert space
	operators,'' in {\em $\uppercase{C}$*-Algebras and Operator Theory} (G.~J.
	Murphy, ed.), pp.~1 -- 76, San Diego: Academic Press, 1990.
	
	\bibitem{ConeAnsari2017}
	Q.~Ansari, E.~K{\"o}bis, and J.~Yao, {\em Vector Variational Inequalities and
		Vector Optimization: Theory and Applications}.
	\newblock Vector Optimization, Springer International Publishing, 2017.
	
	\bibitem{KaplanskyIntroModule}
	I.~Kaplansky, ``Modules over operator algebras,'' {\em American Journal of
		Mathematics}, vol.~75, no.~4, pp.~839--858, 1953.
	
	\bibitem{PaschkeExtModule}
	W.~L. Paschke, ``Inner product modules over \uppercase{$ B $}*-algebras,'' {\em
		Transactions of the American Mathematical Society}, vol.~182, pp.~443--468,
	1973.
	
	\bibitem{RieffelExtModule}
	M.~A. Rieffel, ``Induced representations of $ \uppercase{C} $*-algebras,'' {\em
		Advances in Mathematics}, vol.~13, no.~2, pp.~176 -- 257, 1974.
	
	\bibitem{HilbertModulesChristopher1995}
	E.~C. Lance, ``Chapter 1 - modules and mappings,'' in {\em Hilbert $ C
		$*-Modules: A Toolkit for Operator Algebraists}, London Mathematical Society
	Lecture Note Series (Book 210), pp.~1 -- 13, Cambridge University Press,
	1995.
	
\end{thebibliography}

\end{document}